\documentclass{amsart}
\usepackage{hyperref}

\newcommand{\cB}{\mathcal{B}}
\newcommand{\cC}{\mathcal{C}}
\newcommand{\cG}{\mathcal{G}}
\newcommand{\cJ}{\mathcal{J}}
\newcommand{\cK}{\mathcal{K}}
\newcommand{\cL}{\mathcal{L}}
\newcommand{\cM}{\mathcal{M}}
\newcommand{\cP}{\mathcal{P}}
\newcommand{\cU}{\mathcal{U}}
\newcommand{\cW}{\mathcal{W}}
\newcommand{\cZ}{\mathcal{Z}}
\newcommand{\N}{\mathbb{N}}
\newcommand{\R}{\mathbb{R}}
\newcommand{\Z}{\mathbb{Z}}
\newcommand{\alg}{\mathrm{alg}}
\newcommand{\diag}{\mathrm{diag}}

\newtheorem{theorem}{Theorem}[section]
\newtheorem{lemma}[theorem]{Lemma}
\newtheorem{proposition}[theorem]{Proposition}
\newtheorem{corollary}[theorem]{Corollary}

\begin{document}
\title[Algebras of singular integral operators]
{ALGEBRAS OF SINGULAR INTEGRAL OPERATORS ON NAKANO SPACES
WITH KHVEDELIDZE WEIGHTS OVER CARLESON CURVES
WITH LOGARITHMIC WHIRL POINTS}
\author{A. Yu. Karlovich}
\address{
Universidade do Minho,
Centro de Matem\'atica,
Escola de Ci\^encias,
Campus de Gualtar,
4710-057 Braga,
Portugal}
\email{oleksiy@math.uminho.pt}

\begin{abstract}
We establish a Fredholm criterion for an arbitrary
operator in the Banach algebra of singular integral operators
with piecewise continuous coefficients on Nakano spaces (generalized
Lebesgue spaces with variable exponent) with Khvedelidze weights
over Carleson curves with logarithmic whirl points.
\end{abstract}
\maketitle
\section{Introduction}
Fredholm theory of one-dimensional singular integral operators (SIOs) with
piecewise continuous ($PC$) coefficients on weighted Lebesgue spaces
was constructed by Gohberg and Krupnik \cite{GK92} and \cite{GK70,GK71}
in the beginning of 70s in the case of Khvedelidze weights and piecewise
Lyapunov curves (see also the monographs \cite{CG81,KS01,LS87,MP86}).
Simonenko and Chin Ngok Min \cite{SCNM86} suggested another approach
to the study of the Banach algebra of singular integral operators
with piecewise continuous coefficients on Lebesgue spaces with Khvedelidze
weights over piecewise Lyapunov curves. This approach is based on Simonenko's
local principle \cite{Simonenko65}. In 1992 Spitkovsky \cite{Spitkovsky92}
made a next significant step: he proved a Fredholm criterion for an
individual SIO with $PC$ coefficients on Lebesgue spaces with Muckenhoupt
weights over Lyapunov curves. Finally, B\"ottcher and Yu.~Karlovich extended
Spitkovsky's result to the case of arbitrary Carleson curves and Banach
algebras of SIOs with $PC$ coefficients. After their work the Fredholm
theory of SIOs with $PC$ coefficients is available in the maximal generality
(that, is, when the Cauchy singular integral operator $S$ is bounded on
weighted Lebesgue spaces). We recommend the nice paper \cite{BK01} for a
first reading about this topic and the book \cite{BK97} for a complete
and self-contained exposition.

It is quite natural to consider the same problems in other, more general,
spaces of measurable functions on which the operator $S$ is bounded.
Good candidates for this role are rearrangement-invariant spaces (that is,
spaces with the property that norms of equimeasurable functions are equal).
These spaces have nice interpolation properties and boundedness results
can be extracted from known results for Lebesgue spaces applying interpolation
theorems. The author extended (some parts of) the B\"ottcher-Yu. Karlovich
Fredholm theory of SIOs with $PC$ coefficients to the case of
rearrangement-invariant spaces with Muckenhoupt weights \cite{K98,K02}.
Notice that necessary conditions for the Fredholmness of an individual
singular integral operator with $PC$ coefficients are obtained in \cite{K03}
for weighted reflexive Banach function spaces
on which the operator $S$ is bounded.

Nakano spaces $L^{p(\cdot)}$ (generalized Lebesgue spaces with variable
exponent) are a nontrivial example of Banach function spaces which are not
rearrangement-invariant, in general. Many results about the  behavior of
some classical operators on these spaces have important applications to
fluid dynamics (see \cite{DR03} and the references therein). Kokilashvili
and Samko \cite{KS-GMJ} proved that the operator $S$ is bounded
on weighted Nakano spaces for the case of nice curves, nice weights, and
nice (but variable!) exponents. They also extended the Gohberg-Krupnik Fredholm
criterion for an individual SIO with $PC$ coefficients to this situation
\cite{KS-Proc} (see also \cite{Samko05}). The author \cite{K05} has found a
Fredholm criterion and a formula for the index of an arbitrary operator in
the Banach algebra of SIOs with $PC$ coefficients on Nakano spaces with
Khvedelidze weights over either Lyapunov curves or Radon curves without cusps.

Very recently Kokilashvili and Samko \cite{KS-Memoirs} (see also
\cite[Theorem~7.1]{Kokilashvili05}) have proved a boundedness criterion for
the Cauchy singular integral operator $S$ on Nakano spaces with Khvedelidze
weights over arbitrary Carleson curves. Combining this boundedness result with
the machinery developed in \cite{K03}, we are able to prove a Fredholm criterion
for an individual SIO on a Nakano space with a Khvedelidze weight over a Carleson
curve satisfying a ``logarithmic whirl condition" (see \cite{BK95}, \cite[Ch.~1]{BK97})
at each point. Further, we extend this result to the case of Banach algebras of SIOs
with $PC$ coefficients, using the approach developed in \cite{BK95,BK97,K03,K05}.

The paper is organized as follows. In Section~\ref{sect:preliminaries} we
define weighted Nakano spaces and discuss the boundedness of the operator
$S$ on these spaces. Section~\ref{sect:individual} contains a Fredholm
criterion for an individual SIO with $PC$ coefficients on weighted Nakano
spaces. The proof of this result is based on the local principle of Simonenko
type and factorization technique. In Section~\ref{sect:tools} we formulate
the Allan-Douglas local principle and the two projections theorem. The
results of Section~\ref{sect:tools} are the main tools allowing us to
construct a symbol calculus for the Banach algebra of SIOs with $PC$
coefficients acting on a Nakano space with a Khvedelidze weight over a
Carleson curve with logarithmic whirl points in Section~\ref{sect:symbol}.
\section{Preliminaries}
\label{sect:preliminaries}
\subsection{Weighted Nakano spaces $L_w^{p(\cdot)}$}
Function spaces $L^{p(\cdot)}$ of Lebesgue type with variable exponent
$p$ were studied for the first time by Orlicz  \cite{Orlicz31} in 1931,
but notice that another kind of Banach spaces is called after him.
Inspired by the successful theory of Orlicz spaces, Nakano defined in
the late forties \cite{Nakano50,Nakano51} so-called \textit{modular spaces}.
He considered the space $L^{p(\cdot)}$ as an example of modular spaces.
In 1959, Musielak and Orlicz \cite{MO59} extended the definition of modular
spaces by Nakano. Actually, that paper was the starting point for the theory
of Musielak-Orlicz spaces (generalized Orlicz spaces generated by Young
functions with a parameter), see \cite{Musielak83}.

Let $\Gamma$ be a Jordan (i.e., homeomorphic to a circle) rectifiable
curve. We equip $\Gamma$ with the Lebesgue length measure $|d\tau|$ and
the counter-clockwise orientation. Let $p:\Gamma\to(1,\infty)$ be a measurable
function. Consider the convex modular (see \cite[Ch.~1]{Musielak83} for
definitions and properties)
\[
m(f,p):=\int_\Gamma|f(\tau)|^{p(\tau)}|d\tau|.
\]
Denote by $L^{p(\cdot)}$ the set of all measurable complex-valued functions
$f$ on $\Gamma$ such that $m(\lambda f,p)<\infty$ for some $\lambda=\lambda(f)>0$.
This set becomes a Banach space when equipped with the \textit{Luxemburg-Nakano
norm}
\[
\|f\|_{L^{p(\cdot)}}:=\inf\big\{\lambda>0: \ m(f/\lambda,p)\le 1\big\}
\]
(see, e.g., \cite[Ch.~2]{Musielak83}). Thus, the spaces $L^{p(\cdot)}$ are a
special case of Musielak-Orlicz spaces. Sometimes the spaces $L^{p(\cdot)}$
are referred to as Nakano spaces. We will follow this tradition. Clearly, if
$p(\cdot)=p$ is constant, then the Nakano space $L^{p(\cdot)}$ is isometrically
isomorphic to the Lebesgue space $L^p$. Therefore, sometimes the spaces
$L^{p(\cdot)}$ are called generalized Lebesgue spaces with variable exponent or,
simply, variable $L^p$ spaces.

We shall assume that
\begin{equation}\label{eq:reflexivity}
1<{\rm ess}\inf_{\!\!\!\!\!\!\!\!t\in\Gamma} p(t),
\quad
{\rm ess}\sup_{\!\!\!\!\!\!\!\!\!t\in\Gamma} p(t)<\infty.
\end{equation}
In this case the conjugate exponent
\[
q(t):=\frac{p(t)}{p(t)-1}
\quad (t\in\Gamma)
\]
has the same property.

A nonnegative measurable function $w$ on the curve $\Gamma$ is referred to
as a {\it weight} if $0<w(t)<\infty$ almost everywhere on  $\Gamma$. The
{\it weighted Nakano space} is defined by
\[
L_w^{p(\cdot)}=
\big\{f\mbox{ is measurable on }\Gamma\mbox{ and }fw\in L^{p(\cdot)}\big\}.
\]
The norm in $L_w^{p(\cdot)}$ is defined by
$\|f\|_{L_w^{p(\cdot)}}=\|fw\|_{L^{p(\cdot)}}$.
\subsection{Carleson curves}
A rectifiable Jordan curve $\Gamma$ is said to be a {\it Carleson}
(or {\it Ahlfors-David regular}) {\it curve} if
\[
\sup_{t\in\Gamma}\sup_{R>0}\frac{|\Gamma(t,R)|}{R}<\infty,
\]
where $\Gamma(t,R):=\{\tau\in\Gamma:|\tau-t|<R\}$ for $R>0$
and $|\Omega|$ denotes the measure of a measurable set $\Omega\subset\Gamma$.
We can write
\[
\tau-t=|\tau-t|e^{i\arg(\tau-t)}
\quad\mbox{for}\quad\tau\in\Gamma\setminus\{t\},
\]
and the argument can be chosen so that it is continuous on $\Gamma\setminus\{t\}$.
Seifullaev \cite{Seif80} showed that for an arbitrary Carleson curve the estimate
$\arg(\tau-t)=O(-\log|\tau-t|)$ as $\tau\to t$ holds for every $t\in\Gamma$.
A simpler proof of this result can be found in \cite[Theorem~1.10]{BK97}. One says
that a Carleson curve $\Gamma$ satisfies the \textit{logarithmic whirl condition}
at $t\in\Gamma$ if
\begin{equation}\label{eq:spiralic}
\arg(\tau-t)=-\delta(t)\log|\tau-t|+O(1)\quad (\tau\to t)
\end{equation}
with some $\delta(t)\in\R$. Notice that all piecewise smooth curves satisfy this
condition at each point and, moreover, $\delta(t)\equiv 0$. For more information
along these lines, see \cite{BK95}, \cite[Ch.~1]{BK97}, \cite{BK01}.
\subsection{The Cauchy singular integral operator}
The \textit{Cauchy singular integral} of $f\in L^1$ is defined by
\[
(Sf)(t):=\lim_{R\to 0}\frac{1}{\pi i}\int_{\Gamma\setminus\Gamma(t,R)}
\frac{f(\tau)}{\tau-t}d\tau
\quad (t\in\Gamma).
\]
Not so much is known about the boundedness of the Cauchy singular
integral operator $S$ on weighted Nakano spaces $L_w^{p(\cdot)}$
for general curves, general weights, and general exponents $p(\cdot)$.
{From} \cite[Theorem~6.1]{K03} we immediately get the following.
\begin{theorem}
Let $\Gamma$ be a rectifiable Jordan curve, let $w:\Gamma\to[0,\infty]$
be a weight, and let $p:\Gamma\to(1,\infty)$ be a measurable function
satisfying \eqref{eq:reflexivity}. If the Cauchy singular integral
generates a bounded operator $S$ on the weighted Nakano space
$L_w^{p(\cdot)}$, then
\begin{equation}\label{eq:Ap}
\sup_{t\in\Gamma}\sup_{R>0}\frac{1}{R}
\|w\chi_{\Gamma(t,R)}\|_{L^{p(\cdot)}}
\|\chi_{\Gamma(t,R)}/w\|_{L^{q(\cdot)}}<\infty.
\end{equation}
\end{theorem}

{From} the H\"older inequality for Nakano spaces (see, e.g., \cite{Musielak83} or
\cite{KR91}) and \eqref{eq:Ap} we deduce that if $S$ is bounded on $L_w^{p(\cdot)}$,
then $\Gamma$ is necessarily a Carleson curve. If the exponent
$p(\cdot)=p\in(1,\infty)$ is constant, then \eqref{eq:Ap} is simply the famous
Muckenhoupt condition $A_p$. It is well known that for classical Lebesgue spaces
$L^p$ this condition is not only necessary, but also sufficient for the boundedness
of the Cauchy singular integral operator $S$. A detailed proof of this result can
be found in \cite[Theorem~4.15]{BK97}.

Let $N\in\N$. Consider now a power weight
\begin{equation}\label{eq:power}
\varrho(t):=\prod_{k=1}^N|t-\tau_k|^{\lambda_k},
\quad
\tau_k\in\Gamma,
\quad
k\in\{1,\dots,N\},
\end{equation}
where all $\lambda_k$ are real numbers. Introduce the class $\cP$ of exponents
$p:\Gamma\to(1,\infty)$ satisfying \eqref{eq:reflexivity} and
\begin{equation}\label{eq:Dini-Lipschitz}
|p(\tau)-p(t)|\le\frac{A}{-\log |\tau-t|}
\end{equation}
for some $A\in(0,\infty)$ and all $\tau,t\in\Gamma$ such that $|\tau-t|<1/2$.

The following criterion for the boundedness of the Cauchy singular integral
operator on Nakano spaces with power weights \eqref{eq:Khvedelidze} has been
recently proved by Kokilashvili and Samko \cite{KS-Memoirs} (see also
\cite[Theorem~7.1]{Kokilashvili05}).
\begin{theorem}\label{th:KS}
Let $\Gamma$ be a Carleson Jordan curve, let $\varrho$ be a power weight of
the form \eqref{eq:power}, and let $p\in\cP$. The Cauchy singular integral
operator $S$ is bounded on the weighted Nakano space $L_\varrho^{p(\cdot)}$
if and only if
\begin{equation}\label{eq:Khvedelidze}
0<\frac{1}{p(\tau_k)}+\lambda_k<1
\quad\mbox{for all}\quad
k\in\{1,\dots,N\}.
\end{equation}
\end{theorem}
For weighted Lebesgue spaces over Lyapunov curves the above theorem was proved
by Khvedelidze \cite{Khvedelidze56} (see also the proof in
\cite{GK92,Khvedelidze75,MP86}). Therefore the weights of the form \eqref{eq:power}
are often called \emph{Khvedelidze weights}. We shall follow this tradition.

Notice that if $p$ is constant and $\Gamma$ is a Carleson curve, then
\eqref{eq:Khvedelidze} is equivalent to the fact that $\varrho$ is a Muckenhoupt
weight (see, e.g., \cite[Chapter~2]{BK97}). Analogously one can prove that if the
exponent $p$ belong to the class $\cP$ and the curve $\Gamma$ is Carleson, then
the power weight \eqref{eq:power} satisfies the condition \eqref{eq:Ap} if and only
if \eqref{eq:Khvedelidze} is fulfilled. The proof of this fact is based on certain
estimates for the norms of power functions in Nakano spaces with exponents in the
class $\cP$ (see also \cite[Lemmas~5.7 and 5.8]{K03} and \cite{KS-GMJ}).
\section{Singular integral operators with $PC$ coefficients}
\label{sect:individual}
\subsection{The local principle of Simonenko type}
Let $I$ be the identity operator on $L_\varrho^{p(\cdot)}$.
Under the conditions of Theorem~\ref{th:KS}, the operators
\[
P:=(I+S)/2,
\quad
Q:=(I-S)/2
\]
are bounded projections on $L_\varrho^{p(\cdot)}$ (see \cite[Lemma~6.4]{K03}).
Let $L^\infty$ denote the space of all measurable essentially bounded functions
on $\Gamma$. The operators of the form $aP+Q$ with $a \in L^\infty$ are called
{\it singular integral operators} (SIOs). Two functions $a,b\in L^\infty$ are said
to be locally equivalent at a point $t\in\Gamma$ if
\[
\inf\big\{\|(a-b)c\|_\infty\ :\ c\in C,\ c(t)=1\big\}=0.
\]
\begin{theorem}\label{th:local_principle}
Suppose the conditions of Theorem~{\rm\ref{th:KS}} are satisfied and $a\in L^\infty$.
Suppose for each $t\in\Gamma$ there is a function $a_t\in L^\infty$ which is locally
equivalent to $a$ at $t$. If the operators $a_tP+Q$ are Fredholm on
$L_\varrho^{p(\cdot)}$ for all $t\in\Gamma$, then $aP+Q$ is Fredholm on
$L_\varrho^{p(\cdot)}$.
\end{theorem}
For weighted Lebesgue spaces this theorem is known as Simonenko's local
principle \cite{Simonenko65}. It follows from \cite[Theorem~6.13]{K03}.
\subsection{Simonenko's factorization theorem}
The curve $\Gamma$ divides the complex plane $\mathbb{C}$ into the bounded
simply connected domain $D^+$ and the unbounded domain $D^-$. Without loss
of generality we assume that $0\in D^+$. We say that a function $a\in L^\infty$
admits a \textit{Wiener-Hopf factorization on} $L_\varrho^{p(\cdot)}$ if
$1/a\in L^\infty$ and $a$ can be written in the form
\begin{equation}\label{eq:WH}
a(t)=a_-(t)t^\kappa a_+(t)
\quad\mbox{a.e. on}\ \Gamma,
\end{equation}
where $\kappa\in\Z$, and the factors $a_\pm$ enjoy the following properties:
\begin{enumerate}
\item[{\rm (i)}]
$a_-\in QL_\varrho^{p(\cdot)}\stackrel{\cdot}{+}\mathbb{C}, \quad
1/a_-\in QL_{1/\varrho}^{q(\cdot)}\stackrel{\cdot}{+}\mathbb{C},
\quad a_+\in PL_{1/\varrho}^{q(\cdot)},\quad
1/a_+\in PL_\varrho^{p(\cdot)}$,
\item[{\rm (ii)}]
the operator $(1/a_+)Sa_+I$ is bounded on $L_\varrho^{p(\cdot)}$.
\end{enumerate}
One can prove that the number $\kappa$ is uniquely determined.
\begin{theorem}\label{th:factorization}
Suppose the conditions of Theorem~{\rm\ref{th:KS}} are satisfied. A function
$a\in L^\infty$ admits a Wiener-Hopf factorization {\rm (\ref{eq:WH})} on
$L_\varrho^{p(\cdot)}$ if and only if the operator $aP+Q$ is Fredholm on
$L_\varrho^{p(\cdot)}$. If $aP+Q$ is Fredholm, then its index is equal to
$-\kappa$.
\end{theorem}
This theorem goes back to Simonenko \cite{Simonenko64,Simonenko68}.
For more about this topic we refer to \cite[Section~6.12]{BK97},
\cite[Section~5.5]{BS90}, \cite[Section~8.3]{GK92} and also to \cite{CG81,LS87}
in the case of weighted Lebesgue spaces. Theorem~\ref{th:factorization} follows
from \cite[Theorem~6.14]{K03}.
\subsection{Fredholm criterion for singular integral operators with
$PC$ coefficients}
We denote by $PC$ the Banach algebra of all piecewise continuous functions on
$\Gamma$: a function $a\in L^\infty$ belongs to $PC$ if and only if the finite
one-sided limits
\[
a(t\pm 0):=\lim_{\tau\to t\pm 0}a(\tau)
\]
exist for every $t\in\Gamma$.
\begin{theorem}\label{th:criterion}
Let $\Gamma$ be a Carleson Jordan curve satisfying \eqref{eq:spiralic}
with $\delta(t)\in\R$ for every $t\in\Gamma$. Suppose $p\in\cP$ and $\varrho$
is a power  weight of the form \eqref{eq:power} which satisfies \eqref{eq:Khvedelidze}.
The operator $aP+Q$, where $a\in PC$, is Fredholm on the weighted Nakano space
$L^{p(\cdot)}_\varrho$ if and only if $a(t\pm 0)\ne 0$ and
\begin{equation}
-\frac{1}{2\pi}\arg\frac{a(t-0)}{a(t+0)}
+
\frac{\delta(t)}{2\pi}\log\left|\frac{a(t-0)}{a(t+0)}\right|
+
\frac{1}{p(t)}+\lambda(t)\notin\Z
\label{eq:Fredholm}
\end{equation}
for all $t\in\Gamma$, where
\[
\lambda(t):=\left\{
\begin{array}{lcl}
\lambda_k, &\mbox{if} & t=\tau_k, \quad k\in\{1,\dots,N\},\\
0,         &\mbox{if} & t\notin\Gamma\setminus\{\tau_1,\dots,\tau_N\}.
\end{array}
\right.
\]
\end{theorem}
\begin{proof}
The {\it necessity} part follows from \cite[Theorem~8.1]{K03} because the
B\"ottcher-Yu.~Karlovich indicator functions $\alpha_t$ and $\beta_t$ in that
theorem for a Khvedelidze weight $\varrho$ and a Carleson curve $\Gamma$ satisfying
the logarithmic whirl condition \eqref{eq:spiralic} at $t\in\Gamma$ are calculated by
\[
\alpha_t(x)=\beta_t(x)=\lambda(t)+\delta(t)x\quad \mbox{for}\quad x\in\R
\]
(see \cite[Ch.~3]{BK97} or \cite[Lemma~3.9]{K02}).

{\it Sufficiency.} If $aP+Q$ is Fredholm, then, by \cite[Theorem~6.11]{K03},
$a(t\pm 0)\ne 0$ for all $t\in\Gamma$.

Fix $t\in\Gamma$. For the function $a$ we construct a
``canonical'' function $g_{t,\gamma}$ which is locally equivalent to $a$
at the point $t\in\Gamma$. The interior and the exterior of the unit circle
can be conformally mapped onto $D^+$ and $D^-$ of $\Gamma$, respectively,
so that the point $1$ is mapped to $t$, and the points $0\in D^+$ and
$\infty\in D^-$ remain fixed. Let $\Lambda_0$ and $\Lambda_\infty$
denote the images of $[0,1]$ and $[1,\infty)\cup\{\infty\}$ under this map.
The curve $\Lambda_0\cup\Lambda_\infty$ joins $0$ to $\infty$ and
meets $\Gamma$ at exactly one point, namely $t$. Let $\arg z$ be a
continuous branch of argument in $\mathbb{C}\setminus(\Lambda_0\cup\Lambda_\infty)$.
For $\gamma\in\mathbb{C}$, define the function $z^\gamma:=|z|^\gamma e^{i\gamma\arg z}$,
where $z\in\mathbb{C}\setminus(\Lambda_0\cup\Lambda_\infty)$. Clearly, $z^\gamma$
is an analytic function in $\mathbb{C}\setminus(\Lambda_0\cup\Lambda_\infty)$. The
restriction of $z^\gamma$ to $\Gamma\setminus\{t\}$ will be denoted by
$g_{t,\gamma}$. Obviously, $g_{t,\gamma}$ is continuous and nonzero on
$\Gamma\setminus\{t\}$. Since $a(t\pm 0)\ne 0$, we can define
$\gamma_t=\gamma\in\mathbb{C}$ by the formulas
\[
\operatorname{Re}\gamma_t:=\frac{1}{2\pi}\arg\frac{a(t-0)}{a(t+0)},
\quad
\operatorname{Im}\gamma_t:=-\frac{1}{2\pi}\log\left|\frac{a(t-0)}{a(t+0)}\right|,
\]
where we can take any value of $\arg(a(t-0)/a(t+0))$, which implies that
any two choices of $\operatorname{Re}\gamma_t$ differ by an integer only.
Clearly, there is a constant $c_t\in\mathbb{C}\setminus\{0\}$ such that
$a(t\pm 0)=c_tg_{t,\gamma_t}(t\pm 0)$, which means that $a$ is locally
equivalent to $c_tg_{t,\gamma_t}$ at the point $t\in\Gamma$.
{From} \eqref{eq:Fredholm} it follows that there exists an $m_t\in\Z$ such that
\[
0<m_t-\operatorname{Re}\gamma_t-
\delta(t)\operatorname{Im}\gamma_t+\frac{1}{p(t)}+\lambda(t)<1.
\]
By Theorem~\ref{th:KS}, the operator $S$ is bounded on
$L_{\widetilde{\varrho}}^{p(\cdot)}$, where
\[
\widetilde{\varrho}(\tau)=
|\tau-t|^{m_t-\operatorname{Re}
\gamma_t-\delta(t)\operatorname{Im}\gamma_t}\varrho(\tau)
\]
for $\tau\in\Gamma$. In view of the logarithmic whirl condition \eqref{eq:spiralic}
we have
\begin{eqnarray*}
|(\tau-t)^{m_t-\gamma_t}|
&=&
|\tau-t|^{m_t-\operatorname{Re}\gamma_t}
e^{\operatorname{Im}\gamma_t\arg(\tau-t)}
\\
&=&
|\tau-t|^{m_t-\operatorname{Re}\gamma_t}
e^{-\operatorname{Im}\gamma_t(\delta(t)\log|\tau-t|+O(1))}
\\
&=&
|\tau-t|^{m_t-\operatorname{Re}\gamma_t-\delta(t)\operatorname{Im}\gamma_t}
e^{-\operatorname{Im}\gamma_t O(1)}
\end{eqnarray*}
as $\tau\to t$. Therefore the operator
$\varphi_{t,m_t-\gamma_t}S\varphi_{t,\gamma_t-m_t}I$,
where
\[
\varphi_{t,m_t-\gamma_t}(\tau)=|(\tau-t)^{m_t-\gamma_t}|,
\]
is bounded
on $L_\varrho^{p(\cdot)}$. Then, by \cite[Lemma~7.1]{K03}, the function
$g_{t,\gamma_t}$ admits a Wiener-Hopf factorization on $L_\varrho^{p(\cdot)}$.
Due to Theorem~\ref{th:factorization}, the operator $g_{t,\gamma_t}P+Q$
is Fredholm. Then the operator $c_tg_{t,\gamma_t}P+Q$ is Fredholm, too.
Since the function $c_tg_{t,\gamma_t}$ is locally equivalent to the function
$a$ at every point $t\in\Gamma$, in view of Theorem~\ref{th:local_principle},
the operator $aP+Q$ is Fredholm on $L_\varrho^{p(\cdot)}$.
\end{proof}
\subsection{Double logarithmic spirals}
Given $z_1,z_2\in\mathbb{C}$, $\delta\in\R$, and $r\in(0,1)$, put
\[
\mathcal{S}(z_1,z_2;\delta,r) := \{z_1,z_2\} \cup
\Big\{
z\in\mathbb{C}\setminus\{z_1,z_2\}:
\arg\frac{z-z_1}{z-z_2}-
\delta\log\left|\frac{z-z_1}{z-z_2}\right|\in2\pi (r+\Z)
\Big\}.
\]
The set $\mathcal{S}(z_1,z_2;\delta,r)$ is a double logarithmic spiral whirling
about the points $z_1$ and $z_2$. It degenerates to a familiar Widom-Gohberg-Krupnik
circular arc whenever $\delta=0$ (see \cite{BK97,GK92}).

Fix $t\in\Gamma$ and consider a function $\chi_t\in PC$ which is continuous
on $\Gamma\setminus\{t\}$ and satisfies $\chi_t(t-0)=0$ and $\chi_t(t+0)=1$.

{From} Theorem~\ref{th:criterion} we get the following.
\begin{corollary}\label{co:important}
Let $\Gamma$ be a Carleson Jordan curve satisfying \eqref{eq:spiralic}
with $\delta(t)\in\R$ for every $t\in\Gamma$. Suppose $p\in\cP$ and $\varrho$ is
a power weight of the form \eqref{eq:power} which satisfies \eqref{eq:Khvedelidze}.
Then
\[
\big\{\lambda\in\mathbb{C}:(\chi_t-\lambda)P+Q
\mbox{ is not Fredholm on }L_\varrho^{p(\cdot)}\big\}
=\mathcal{S}\big(0,1;\delta(t),1/p(t)+\lambda(t)\big).
\]
\end{corollary}
\section{Tools for the construction of the symbol calculus}\label{sect:tools}
\subsection{The Allan-Douglas local principle}
Let $B$ be a Banach algebra with identity. A subalgebra $Z$ of $B$ is said to
be a central subalgebra if $zb=bz$ for all $z\in Z$ and all $b\in B$.
\begin{theorem}\label{th:AllanDouglas}
{\rm (see \cite[Theorem~1.34(a)]{BS90}).}
Let $B$ be a Banach algebra with unit $e$ and let $Z$ be closed central subalgebra
of $B$ containing $e$. Let $M(Z)$ be the maximal ideal space of $Z$, and for
$\omega\in M(Z)$, let $J_\omega$ refer to the smallest closed two-sided ideal
of $B$ containing the ideal $\omega$. Then an element $b$ is invertible in $B$
if and only if $b+J_\omega$ is invertible in the quotient algebra $B/J_\omega$
for all $\omega\in M(Z)$.
\end{theorem}
\subsection{The two projections theorem}
The following two projections theorem was obtained by Finck,
Roch, Silbermann \cite{FRS93} and Gohberg, Krupnik \cite{GK93}.
\begin{theorem}\label{th:2proj}
Let $F$ be a Banach algebra with identity $e$, let
$\cC$ be a Banach subalgebra of $F$ which
contains $e$ and is isomorphic to ${\mathbb{C}}^{n \times n}$,
and let $p$ and $q$ be two projections in $F$ such that
$cp=pc$ and $cq=qc$ for all $c \in \cC$. Let $W=\alg(\cC,p,q)$
be the smallest closed subalgebra of $F$ containing $\cC,p,q$. Put
\[
x=pqp+(e-p)(e-q)(e-p),
\]
denote by $\mathrm{sp}\,x$ the spectrum of $x$ in $F$, and suppose the points
$0$ and $1$ are not isolated points of $\mathrm{sp}\,x$. Then
\begin{enumerate}
\item[{\rm (a)}]
for each $\mu \in \mathrm{sp}\,x$ the map $\sigma_{\mu}$ of $\cC \cup \{p,q\}$
into the algebra ${\mathbb{C}}^{2n\times 2n}$ of all complex $2n\times 2n$
matrices defined by
\begin{equation}\label{eq:2proj1}
\sigma_{\mu}c=\left(
\begin{array}{cc}
c & 0\\
0 & c
\end{array}
\right),
\quad
\sigma_{\mu}p=\left(
\begin{array}{cc}
E & 0\\
0 & 0
\end{array}
\right),
\end{equation}
\begin{equation}\label{eq:2proj2}
\sigma_{\mu}q=\left(
\begin{array}{cc}
\mu E & \sqrt{\mu(1-\mu)}E \\
\sqrt{\mu(1-\mu)}E  & (1-\mu)E
\end{array}
\right),
\end{equation}
where $c\in \cC, E$ denotes the $n \times n$ unit matrix and $\sqrt{\mu(1-\mu)}$
denotes any complex number whose square is $\mu(1-\mu)$, extends to a Banach algebra
homomorphism
\[
\sigma_{\mu}: W \to {\mathbb{C}}^{2n \times 2n};
\]

\item[{\rm (b)}]
every element $a$ of the algebra $W$ is invertible in the algebra $F$ if and only if
\[
\det \sigma_{\mu} a \neq 0 \quad\mbox{for all}\quad \mu \in \mathrm{sp}\,x;
\]

\item[{\rm (c)}]
the algebra $W$ is inverse closed in $F$ if and only if the spectrum of
$x$ in $W$ coincides with the spectrum of $x$ in $F$.
\end{enumerate}
\end{theorem}

A further generalization of the above result to the case of $N$ projections
is contained in \cite{BK97}.
\section{Algebra of singular integral operators with $PC$ coefficients}
\label{sect:symbol}
\subsection{The ideal of compact operators}
In this section we will suppose that $\Gamma$ is a Carleson curve satisfying
\eqref{eq:spiralic} with $\delta(t)\in\R$ for every $t\in\Gamma$, $p\in\cP$,
and $\varrho$ is a Khvedelidze weight of the form \eqref{eq:power} which
satisfies \eqref{eq:Khvedelidze}.
Let $X_n:=[L_\varrho^{p(\cdot)}]_n$ be the direct sum of $n$ copies of weighted
Nakano spaces $X:=L_\varrho^{p(\cdot)}$, let $\cB:=\cB(X_n)$ be the Banach algebra
of all bounded linear operators on $X_n$, and let $\cK:=\cK(X_n)$ be the closed
two-sided ideal of all compact operators on $X_n$. We denote by $C^{n\times n}$
(resp. $PC^{n\times n}$) the collection of all continuous (resp. piecewise
continuous) $n\times n$ matrix functions, that is, matrix-valued functions with
entries in $C$ (resp. $PC$). Put $I^{(n)}:=\diag\{I,\dots, I\}$ and
$S^{(n)}:=\diag\{S,\dots,S\}$. Our aim is to get a Fredholm criterion for an operator
\[
A\in\cU:=\alg(PC^{n\times n},S^{(n)}),
\]
the smallest Banach subalgebra of $\cB$
which contains all operators of multiplication by matrix-valued functions in
$PC^{n\times n}$ and the operator $S^{(n)}$.
\begin{lemma}\label{le:compact}
The ideal
$\cK$ is contained in the algebra $\alg(C^{n\times n},S^{(n)})$, the smallest closed
subalgebra of $\cB$ which contains the operators of multiplication by
continuous matrix-valued functions and the operator $S^{(n)}$.
\end{lemma}

The proof of this statement is standard and can be developed as in
\cite[Lemma~9.1]{K96} or \cite[Lemma~5.1]{K05}.
\subsection{Operators of local type}
We shall denote by $\cB^\pi$ the Calkin algebra $\cB/\cK$ and by $A^\pi$
the coset $A+\cK$ for any operator $A\in\cB$. An operator $A\in\cB$ is said
to be of {\it local type} if $AcI^{(n)}-cA$ is compact for all $c\in C$,
where $cI^{(n)}$ denotes the operator of multiplication by the diagonal
matrix-valued function $\diag\{c,\dots,c\}$. This notion goes back to Simonenko
\cite{Simonenko65} (see also the presentation of Simonenko's local theory
in his joint monograph with Chin Ngok Min \cite{SCNM86}). It easy to see that
the set $\cL$ of all operators of local type is a closed subalgebra of $\cB$.
\begin{proposition}\label{pr:OLT}
\begin{enumerate}
\item[{\rm (a)}] We have
$\cK\subset\cU\subset\cL$.

\item[{\rm (b)}]
An operator $A\in\cL$ is Fredholm if and only if the coset $A^\pi$ is invertible
in the quotient algebra $\cL^\pi:=\cL/\cK$.
\end{enumerate}
\end{proposition}
\begin{proof}
(a) The embedding $\cK\subset\cU$ follows from Lemma~\ref{le:compact},
the embedding $\cU\subset\cL$ follows from the fact that $cS-ScI$ is
a compact operator on $L_\varrho^{p(\cdot)}$ for $c\in C$
(see, e.g., \cite[Lemma~6.5]{K03}).

\medskip
(b) The proof of this fact is straightforward.
\end{proof}
\subsection{Localization}
{From} Proposition~\ref{pr:OLT}(a) we deduce that the quotient algebras
$\cU^\pi:=\cU/\cK$ and $\cL^\pi:=\cL/\cK$ are well defined. We shall
study the invertibility of an element $A^\pi$ of $\cU^\pi$ in the larger
algebra $\cL^\pi$ by using a localization techniques (more precisely,
Theorem~\ref{th:AllanDouglas}). To this end, consider
\[
\cZ^\pi:=\big\{(cI^{(n)})^\pi:c\in C\big\}.
\]
{From} the definition of $\cL$ it follows that $\cZ^\pi$ is a central subalgebra
of $\cL^\pi$. The maximal ideal space $M(\cZ^\pi)$ of $\cZ^\pi$ may be
identified with the curve $\Gamma$ via the Gelfand map $\cG$ given by
\[
\cG:\cZ^\pi\to C,
\quad
\big(\cG(cI^{(n)})^\pi\big)(t)=c(t)
\quad (t\in\Gamma).
\]
In accordance with Theorem~\ref{th:AllanDouglas}, for every $t\in\Gamma$
we define $\cJ_t\subset\cL^\pi$ as the smallest closed two-sided ideal of
$\cL^\pi$ containing the set
\[
\big\{(cI^{(n)})^\pi\ :\ c\in C,\ c(t)=0\big\}.
\]

Consider a function $\chi_t\in PC$ which is continuous
on $\Gamma\setminus\{t\}$ and satisfies $\chi_t(t-0)=0$ and $\chi_t(t+0)=1$.
For $a\in PC^{n\times n}$ define the function $a_t\in PC^{n\times n}$
by
\begin{equation}\label{eq:at}
a_t:=a(t-0)(1-\chi_t)+a(t+0)\chi_t.
\end{equation}
Clearly $(aI^{(n)})^\pi-(a_tI^{(n)})^\pi\in\cJ_t$. Hence, for any operator
$A\in\cU$, the coset $A^\pi+\cJ_t$ belongs to the smallest closed subalgebra
$\cW_t$ of $\cL^\pi/\cJ_t$ containing the cosets
\begin{equation}\label{eq:projections}
p:=\big((I^{(n)}+S^{(n)})/2\big)^\pi+\cJ_t,
\
q:=(\chi_tI^{(n)})^\pi+\cJ_t,
\end{equation}
where $\chi_tI^{(n)}$ denotes the operator of multiplication by the diagonal
matrix-valued function $\diag\{\chi_t,\dots,\chi_t\}$ and the algebra
\begin{equation}\label{eq:algebra}
\cC:=\big\{(cI^{(n)})^\pi+\cJ_t\ : \ c\in\mathbb{C}^{n\times n}\big\}.
\end{equation}
The latter algebra is obviously isomorphic to $\mathbb{C}^{n\times n}$,
so $\cC$ and $\mathbb{C}^{n\times n}$ can be identified with each other.
\subsection{The spectrum of $pqp+(e-p)(e-q)(e-p)$}
Since $P^2=P$ on $L_\varrho^{p(\cdot)}$ (see, e.g., \cite[Lemma~6.4]{K03})
and $\chi_t^2-\chi_t\in C$, $(\chi_t^2-\chi_t)(t)=0$,
it is easy to see that
\begin{equation}\label{eq:2proj-conditions}
p^2=p,
\quad
q^2=q,
\quad
pc=cp,
\quad
qc=cq
\end{equation}
for every $c\in\cC$, where $p,q$ and $\cC$ are given by \eqref{eq:projections}
and \eqref{eq:algebra}. To apply Theorem~\ref{th:2proj} to the algebras
$F=\cL^\pi/\cJ_t$ and $W=\cW_t=\alg(\cC,p,q)$, we have to identify the spectrum
of
\begin{equation}
pqp+(e-p)(e-q)(e-p)
=\big(P^{(n)}\chi_tP^{(n)}+Q^{(n)}(1-\chi_t)Q^{(n)}\big)^\pi+\cJ_t
\label{eq:element}
\end{equation}
in the algebra $F=\cL^\pi/\cJ_\tau$; here $P^{(n)}:=(I^{(n)}+S^{(n)})/2$
and $Q^{(n)}:=(I^{(n)}-S^{(n)})/2$.
\begin{lemma}\label{le:spectrum}
Let $\chi_t\in PC$ be a continuous function on $\Gamma\setminus\{t\}$ such that
\[
\chi_t(t-0)=0,
\quad
\chi_t(\tau+0)=1,
\]
and
\[
\chi_t(\Gamma\setminus\{t\})\cap\mathcal{S}(0,1;\delta(t),1/p(t)+\lambda(t))=\emptyset.
\]
Then the spectrum of \eqref{eq:element} in $\cL^\pi/\cJ_t$ coincides
with $\mathcal{S}(0,1;\delta(t),1/p(t)+\lambda(t))$.
\end{lemma}
\begin{proof}
Once we have Corollary~\ref{co:important}  at hand, the proof of this lemma
can be developed by a literal repetition of the proof of \cite[Lemma~9.4]{K96}.
\end{proof}
\subsection{Symbol calculus}
Now we are in a position to prove the main result of this paper.
\begin{theorem}\label{th:symbol}
Define the ``double logarithmic spirals bundle''
\[
\cM:=
\bigcup\limits_{t\in\Gamma} \Big(\{t\} \times
\mathcal{S}\big(0,1;\delta(t),1/p(t)+\lambda(t)\big) \Big).
\]
\begin{enumerate}
\item[{\rm (a)}]
For each point $(t,\mu)\in\cM$, the map
\[
\sigma_{t,\mu} \: : \:
\{S^{(n)}\}\cup\{aI^{(n)}\: :\:
a\in PC^{n\times n}\} \to \mathbb{C}^{2n\times 2n}
\]
given by
\[
\sigma_{t,\mu}(S^{(n)})
=
\left(
\begin{array}{ll}
E &  O\\
O & -E
\end{array}
\right),
\
\sigma_{t,\mu}(aI^{(n)})
=
\left(
\begin{array}{ll}
a_{11}(t,\mu)  &  a_{12}(t,\mu)\\
a_{21}(t,\mu)  &  a_{22}(t,\mu)
\end{array}
\right),
\]
where
\begin{eqnarray*}
a_{11}(t,\mu)
&:=&
a(t+0)\mu + a(t-0)(1-\mu),\\
a_{12}(t,\mu)
&=&
a_{21}(t,\mu)
:=
(a(t+0)-a(t-0)) \sqrt{\mu(1-\mu)},
\\
a_{22}(t,\mu)
&:=&
a(t+0)(1-\mu) + a(t-0)\mu,
\end{eqnarray*}
and $O$ and $E$ are the zero and identity $n\times n$ matrices, respectively,
extends to a Banach algebra homomorphism
\[
\sigma_{t,\mu} :\cU\to\mathbb{C}^{2n\times 2n}
\]
with the property that $\sigma_{t,\mu}(K)$ is the zero matrix
for every compact operator $K$ on $X_n$;

\item[{\rm (b)}]
an operator $A\in\cU$ is Fredholm on $X_n$
if and only if
\[
\det\sigma_{t,\mu} (A)\neq 0
\quad\mbox{for all}\quad (t,\mu)\in\cM;
\]

\item[{\rm (c)}]
the quotient algebra $\cU^\pi$ is inverse closed in the
Calkin algebra $\cB^\pi$, that is, if a coset $A^\pi\in\cU^\pi$
is invertible in $\cB^\pi$, then $(A^\pi)^{-1}\in\cU^\pi$.
\end{enumerate}
\end{theorem}
\begin{proof}
The idea of the proof of this theorem is borrowed from \cite{BK97}
and is based on the Allan-Douglas local principle
(Theorem~\ref{th:AllanDouglas}) and the two projections theorem
(Theorem~\ref{th:2proj}).

Fix $t\in\Gamma$ and choose a function $\chi_t\in PC$ such that
$\chi_t$ is continuous on $\Gamma\setminus\{t\}$, $\chi_t(t-0)=0$,
$\chi_t(t+0)=1$, and
$\chi_t(\Gamma\setminus\{t\})\cap\mathcal{S}(0,1;\delta(t),1/p(t)+\lambda(t))=\emptyset$.
{From} \eqref{eq:2proj-conditions} and Lemma~\ref{le:spectrum} we deduce that
the algebras $\cL^\pi/\cJ_t$ and $\cW_t=\alg(\cC,p,q)$, where $p,q$ and $\cC$
are given by \eqref{eq:projections} and \eqref{eq:algebra}, respectively,
satisfy all the conditions of the two projections theorem (Theorem~\ref{th:2proj}).

(a) In view of Theorem~\ref{th:2proj}(a), for every
$\mu\in\mathcal{S}(0,1;\delta(t),1/p(t)+\lambda(t))$,
the map $\sigma_\mu:\mathbb{C}^{n\times n}\cup\{p,q\}\to\mathbb{C}^{2n\times 2n}$
given by \eqref{eq:2proj1}--\eqref{eq:2proj2} extends to a Banach algebra homomorphism
$\sigma_\mu:\cW_t\to\mathbb{C}^{2n\times 2n}$. Then the map
\[
\sigma_{t,\mu}=\sigma_\mu\circ\pi_t:\cU\to\mathbb{C}^{2n\times 2n},
\]
where $\pi_t:\cU\to\cW_t=\cU^\pi/\cJ_t$ is acting by the rule $A\mapsto A^\pi+\cJ_t$,
is a well defined Banach algebra homomorphism and
\[
\sigma_{t,\mu}(S^{(n)})=2\sigma_\mu p-\sigma_\mu e=
\left(\begin{array}{cc}E & O \\O & -E\end{array}\right).
\]
If $a\in PC^{n\times n}$, then in view of \eqref{eq:at} and
$(aI^{(n)})^\pi-(a_tI^{(n)})^\pi\in\cJ_t$ it follows that
\begin{eqnarray*}
\sigma_{t,\mu}(aI^{(n)})
&=&
\sigma_{t,\mu}(a_tI^{(n)})
=
\sigma_\mu(a(t-0))\sigma_\mu(e-q)+\sigma_\mu(a(t+0))\sigma_\mu q
\\
&=&
\left(\begin{array}{cc} a_{11}(t,\mu) & a_{12}(t,\mu)\\
a_{21}(t,\mu) & a_{22}(t,\mu)\end{array}\right).
\end{eqnarray*}
{From} Proposition~\ref{pr:OLT}(a) it follows that $\pi_t(K)=K^\pi+\cJ_t=\cJ_t$ for every
$K\in\cK$ and every $t\in\Gamma$. Hence,
\[
\sigma_{t,\mu}(K)=\sigma_\mu(0)=
\left(\begin{array}{cc} O & O\\ O & O\end{array}\right).
\]
Part (a) is proved.

(b) {From} Proposition~\ref{pr:OLT} it follows that the Fredholmness of $A\in\cU$ is
equivalent to the invertibility of $A^\pi\in\cL^\pi$. By Theorem~\ref{th:AllanDouglas},
the former is equivalent to the invertibility of $\pi_t(A)=A^\pi+\cJ_t$
in $\cL^\pi/\cJ_t$ for every $t\in\Gamma$. By Theorem~\ref{th:2proj}(b),
this is equivalent to
\begin{eqnarray}\label{eq:symbol1}
\det\sigma_{t,\mu}(A)=\det\sigma_\mu\pi_t(A)\ne 0
\mbox{ for all }(t,\mu)\in\cM.
&&
\end{eqnarray}
Part (b) is proved.

(c) Since $\mathcal{S}(0,1;\delta(t),1/p(t)+\lambda(t))$ does not separate
the complex plane $\mathbb{C}$, it follows that the spectra of
\eqref{eq:element} in the algebras
$\cL^\pi/\cJ_t$ and $\cW_t=\cU^\pi/\cJ_t$ coincide, so we can apply
Theorem~\ref{th:2proj}(c). If $A^\pi$, where $A\in\cU$, is invertible in $\cB^\pi$,
then \eqref{eq:symbol1} holds. Consequently, by Theorem~\ref{th:2proj}(b), (c),
$\pi_t(A)=A^\pi+\cJ_t$ is invertible in $\cW_t=\cU^\pi/\cJ_t$ for every $t\in\Gamma$.
Applying Theorem~\ref{th:AllanDouglas} to $\cU^\pi$, its central subalgebra $\cZ^\pi$,
and the ideals $\cJ_t$, we obtain that $A^\pi$ is invertible in $\cU^\pi$, that is,
$\cU^\pi$ is inverse closed in the Calkin algebra $\cB^\pi$.
\end{proof}

Note that the approach to the study of Banach algebras of SIOs
based on the Allan-Douglas local principle and the two projections
theorem is nowadays standard. It was successfully applied in many
situations (see, e.g., \cite{BK95,BK97,FRS93,K96,K98,K02,K05}).
However, it does not allow to get formulas for the index of an arbitrary
operator in the Banach algebra of SIOs with $PC$ coefficients.
These formulas can be obtained similarly to the classical situation
considered by Gohberg and Krupnik \cite{GK70,GK71} (see also \cite[Ch.~10]{BK97}).
For reflexive Orlicz spaces over Carleson curves with logarithmic whirl points
this was done by the author \cite{K98-index}. In the case of Nakano spaces with
Khvedelidze weights over Carleson curves with logarithmic whirl points the
index formulas are almost the same as in \cite{K98-index}. It is only necessary to
replace the both Boyd indices $\alpha_M$ and $\beta_M$ of an Orlicz space $L^M$
by the numbers $1/p(t)+\lambda(t)$ in corresponding formulas.

\newcommand{\ieot}{Integral Equations Operator Theory}
\newcommand{\otaa}{Operator Theory: Advances and Applications}

\end{document}